\documentclass[a4paper]{article} 

\usepackage{graphicx}
\usepackage{epstopdf}
\usepackage{amsfonts}
\usepackage{amsmath}
\usepackage{numprint}
\npdecimalsign{.}

\def\RR{{\mathbb R}}

\begin{document}

\title{Accurate Spectral Collocation Computation of High Order Eigenvalues for Singular Schr\"odinger Equations}

\author{C\u{a}lin-Ioan Gheorghiu
\thanks{Tiberiu Popoviciu Institute of Numerical Analysis, Romanian Academy; ghcalin@ictp.acad.ro}}

\date{ }

\maketitle

\abstract{We are concerned with the study of some classical spectral collocation methods as well as with the new software system Chebfun in computing high order eigenpairs of singular and regular Schr\"odinger eigenproblems.
We want to highlight both the qualities as well as the shortcomings of these methods and evaluate them in conjunction with the usual ones.
In order to resolve a boundary singularity we use Chebfun with domain truncation. Although it is applicable with spectral collocation, a special technique to introduce boundary conditions as well as a coordinate transform, which maps an unbounded domain to a finite one, are the special ingredients.
A challenging set of \textquotedblleft{hard}\textquotedblright benchmark problems, for which usual numerical methods (f. d., f. e. m., shooting etc.) fail, are analyzed. In order to separate \textquotedblleft good\textquotedblright and \textquotedblleft bad\textquotedblright eigenvalues we estimate the drift of the set of eigenvalues of interest with respect to the order of approximation and/or scaling of domain parameter. It automatically provides us with a measure of the error within which the eigenvalues are computed and a hint on numerical stability.
We pay a particular attention to problems with almost multiple eigenvalues as well as to problems with a mixed  spectrum.


}

\section{Introduction}
There is clearly an increasing interest to develop accurate and efficient methods of solution to singular Schr\"odinger eigenproblems. Our main interest here is to compare the capabilities of the new Chebfun package with those of classical spectral methods
in solving such problems, showing their capabilities and weaknesses. The latter
employ basis functions and/or grid points based on Chebyshev, Laguerre or Hermite polynomials as well as on sinc or Fourier functions. The effort expended by both classes of methods is also of real interest. It can be assessed in terms of the ease of implementation of the methods as well as in terms of computer resources required to achieve a specified accuracy.

Spectral methods have been shown to provide exponential convergence for a large variety of problems, generally with smooth solutions, and are often preferred. For details on Chebfun we refer to \cite{DBT}, \cite{DHT}, \cite{DHT1}, \cite{OT}, \cite{TBD} and \cite{LNT}.  For Chebyshev collocation (ChC),  Laguerre-Gauss-Radau collocation (LGRC), Hermite (HC) and sinc collocation (SiC) we refer among other sources to our contributions \cite{Cig14}, \cite{Cig18} as well as to the seminal paper \cite{WR}.

For problems on the entire real axis SiC proved to be particularly well suited. Moreover, this method has given excellent results recorded in our contribution \cite{Cig18} and in the works cited there. The so-called generalized pseudospectral (GPS) method, actually the Legendre collocation, is employed in \cite{Roy} to calculate the
bound states of the Hulth\'{e}n and the Yukawa potentials in quantum mechanics, with special
emphasis on higher excited states and stronger couplings. The author uses a two parameter dependent nonlinear transformation in order to map the half-line into the canonical interval $\left[-1, 1\right].$ In contrast, we will use an analogous transformation but which depends on only one parameter. Also, very recently spectral methods based on non-classical orthogonal polynomials have been used in \cite{BDS} in order to solve some Schr\"odinger problems connected with Fokker-Planck operator.
A particular attention will be paid in this paper to the challenging issue of \textit{continuous spectra vs. discrete
(numerical) eigenvalues}. It is well known that some eigenvalue problems (see for instance the well known text \cite{BR}) for differential operators which are
naturally posed on the whole real line or the half-line,
often lead to some discrete eigenvalues plus a continuous spectrum. Actually, the usual numerical approximation
typically involves three processes:
\begin{enumerate}
\item reduction to a finite interval;
\item discretization;
\item application of a numerical eigenvalue solver.
\end{enumerate}
Reduction to a finite interval and discretization typically eliminate the continuous
spectrum. Even if we do not truncate a priory the domain on which the problem is
formulated such an inherent reduction can not be avoided.

It can be argued that, generally speaking, in solving various differential problems, the Chebfun software provides a greater flexibility than the classical spectral methods. This fact is fully true for regular problems. 

Unfortunately, in the presence of various singularities, the maximum order of approximation $N,$ of the unknowns can be reached ($N \geq 4000$) and then Chebfun issues a message that warns about the possible inaccuracy of the results provided.

We came out of this tangle using modified classical spectral methods. In this way, when we had serious doubts about the accuracy of the solutions given by Chebfun, we managed to establish the correctness of the numerical results.

As a matter of fact, in order to resolve a singularity on the ends of an unbounded integration interval, Chebfun uses only the arbitrary truncation of the domain. Classical spectral methods can also use this method, but it is not recommended. For singular points at finite distances (mainly origin) we will use the so-called \textit{removing technique of independent boundary conditions.} The boundary conditions at infinity can be enforced using basis functions that satisfy these conditions (Laguerre, Hermite, sinc). An alternative method, which proved to be very accurate, is the Chebyshev collocation (ChC) method in combination with a change of
variables (coordinates) which transform the half line into the canonical Chebyshev interval $\left[-1, 1\right].$ Then the removing technique of independent boundary conditions is essential in order to remove the singularities at the end points of integration interval.

A Chebfun code and two MATLAB codes, one for ChC and another for SiC method, are provided in order to exemplify. With minor modifications they could be fairly useful for various numerical experiments.

The structure of this work is as follows. In Section \ref{SP} we recall some specific issues for the regular as well as singular Schr\"odinger eigenproblems. The main comment refers to the notion of \textit{mixed spectrum}. In Section \ref{SC} we review on the Chebfun structure and the classical spectral methods (differentiation matrices, enforcing boundary conditions, etc.). The Section \ref{numerics} is the central part of the paper. Here we analyze some benchmark problems. In order to separate the \textquotedblleft{good}\textquotedblright from the \textquotedblleft{bad}\textquotedblright eigenvalues we estimate their relative drift with respect to some parameters. The accuracy in computing eigenfunctions is estimated by their departure from orthogonality. We end up with Section \ref{conclusions} where we underline some conclusions and suggest some open problems.
\section{Regular and singular Schr\"odinger eigenproblems} \label{SP}
The Schr\"{o}dinger equation reads%
\begin{equation}
u^{\prime \prime }+\left[ \lambda -q\left( x\right) \right] u=0,\ -\infty
\leq a<x<b\leq +\infty .  \label{S_eq}
\end{equation}%
It is a Liouville normal form of a general Sturm-Liouville (SL) equation where $%
\lambda $ is proportional with the \textit{energy levels} of the physical
system, $q\left( x\right) $ is directly proportional with the
\textit{potential energy} and the \textquotedblleft wave
function\textquotedblright $u$ may be real or complex such that $u\ u^{\ast
}dx=\left\vert u\right\vert ^{2}dx$ is the probability that the particle
under consideration will be \textquotedblleft observed\textquotedblright in
the interval $\left( x,\ x+dx\right) .$ In problems involving Schr\"{o}%
dinger equations, it is customary among chemists and physicists to define
the \textit{spectrum }of this Sturm-Liouville problem as the all eigenvalues
$\lambda $ for which eigenfunctions $u$ exist. The set of isolated points
(if any) in this spectrum is called the \textit{discrete spectrum; }the part
(if any) that consists of entire interval is called \textit{continuous
spectrum.} We shall adopt this suggestive terminology here.

The equation (\ref{S_eq}) can be given on a finite, semi-infinite, or
infinite interval. Only on a closed and finite interval $a\leq x\leq b$ can
Schr\"{o}dinger equation be associated with a \textit{regular }%
Sturm-Liouville problem. If the interval of definition is semi-infinite or
infinite, or is finite and $q\left( x\right) $ vanishes at one or both
endpoints, or if $q$ \ is discontinuous, we can not obtain from  (\ref{S_eq}%
) a regular Sturm-Liouville problem. In any such case, the Schr\"{o}dinger
equation (\ref{S_eq}) is called \textit{singular.} We obtain a \textit{singular} \textit{%
eigenproblem} from a singular Schr\"{o}dinger equation by imposing suitable
homogeneous boundary conditions. They can not always be described by
formulae like $\alpha u\left( e\right) +\beta u^{\prime }\left( e\right) =0,$
where $e$ can be the end point $a$ or $b.$ For instance, the condition that $%
u$ be bounded near a singular end point, which can be finite or $\pm \infty,$ is a common boundary condition defining a singular eigenproblem.
For regular SL problems, it is proved (see for instance \cite%
{BR}) that the spectrum is always discrete, and the eigenfunctions are
(trivially) square-integrable.
For singular problems the situation is completely different.
For instance in their textbook \cite{BR} Birkhoff and Rota consider the eigenproblem
attached to the \textit{free particle} equation
and show that the spectrum of the free particle is continuous.
Some software packages have been designed over time to solve various singular SL problems. The most important would be
SLEIGN and SLEIGN2, SLEDGE, SL02F and MATSLISE. The SLDRIVER interactive package
supports exploration of a set of SL problems with the four previously mentioned packages.
In \cite{PF} (see also \cite{PFX}) the authors designed the software package SLEDGE. They  observed that for a class of
singular problems their method either fails or converges very slowly. Essentially, the numerical
method used in this software package replaces the coefficient function $q(x)$ by step function approximation. Similar behavior has been observed on the NAG code SL02F introduced in \cite{MP} and \cite{PM} as well as on the packages SLEIGN
and SLEIGN2 introduced in \cite{BEZ} and \cite{BGKZ}. The  MATSLISE code introduced in  \cite{Ledoux} can solve some
Schr\"{o}dinger eigenvalue problem by a constant perturbation method of a higher order.

The main purpose of this paper is to argue that Chebfun, along with the spectral collocation methods, can be a very feasible alternative to these software packages regarding accuracy, robustness as well as simplicity of implementation. In addition, these methods can compute the \textquotedblleft whole\textquotedblright set of eigenvectors and provide some details on the accuracy and numerical stability of the results provided.

\section{Chebfun vs. spectral collocation (ChC, LGRC, SiC)} \label{SC}

\subsection{Chebfun}
The Chebfun system, in object-oriented MATLAB, contains algorithms which amount to
spectral collocation methods on Chebyshev grids of automatically determined resolution.  Its properties are briefly summarized in \cite{DHT}. In \cite{DBT} the authors explain that \textit{chebops} are the fundamental Chebfun tools for solving ordinary differential (or integral) equations. One may then use them as
tools for more complicated computations that may be nonlinear and may involve
partial differential equations. This is analogous to the situation in MATLAB itself.
The implementation of chebops combines the numerical analysis idea of spectral
collocation with the computer science idea of \textit{lazy or delayed evaluation
of the associated spectral discretization matrices}. The grammar of chebops along with a lot of illustrative  examples is displayed in the above quoted paper as well as in the text \cite{TBD}. Thus one can get a suggestive image of what they can do.

In \cite{DBT} p.12 the authors explain clearly how the Chebfun works, i.e., it solves the eigenproblem for two different orders of approximation, automatically chooses a reference eigenvalue  and checks the convergence of the process. At the same time, it warns about the possible failures due to the high non-normality of the analyzed operator (matrix).

Actually, we want to show in this paper that Chebfun along with chebops can do much more, i.e., can accurately solve highly (double) singular Schr\"{o}dinger eigenproblems.

\subsection{ChC, LGRC and SiC}
In the spectral collocation method the unknown solution to a differential equation is expanded as a global interpolant, such as a trigonometric or polynomial interpolant. In other methods, such as f. e. and f. d., the underlying expansion involves local interpolants such as piecewise polynomials. This means that the accuracy of spectral collocation is superior. For problems with smooth solutions convergence rates are typically of order $e^{-cN}$ or $e^{-c\sqrt{N}}$ where $N$ is the order of approximation or resolution, i.e., the number of degree of freedom in expansion. In contrast, f. e. or f. d. yield convergence rates that are only algebraic in $N$, typically of orders $N^{-2}$ or $N^{-4}.$
The net superiority of global spectral methods on local methods is discussed in detail in \cite{ST}.
In all spectral collocation methods designed so far we have used the collocation differentiation matrices from the seminal paper \cite{WR}. We preferred this MATLAB differentiation suite for the accuracy, efficiency as well as for the ingenious way of introducing various boundary conditions.

In order to impose (enforce) the boundary conditions we have used two methods that are conceptually different, namely the \textit{boundary bordering} as well as the \textit{basis recombination}. A very efficient way to accomplish the boundary bordering is available in \cite{JH} and is called \textit{removing technique of independent boundary conditions}. We have used this technique in the large majority of our papers except \cite{GigIsp} where the latter technique has been employed. In the last quoted paper a modified Chebyshev tau method based on basis recombination has been used in order to solve an Orr-Sommerfeld problem with an eigenparameter dependent boundary condition. Even in eigenproblems that contain the spectral parameter in the definition of the boundary conditions we have used the boundary bordering technique (see \cite{Cig17}).

In \cite{CIGHPR} (see also \cite{PGH}) we have solved some multiparameter eigenproblems (MEP) which come from separation of variables, in several orthogonal coordinate systems, applied to the Helmholtz, Laplace, or Schr\"{o}dinger equation. Important cases include Mathieu's system, Lame's system, and a system of spheroidal wave functions.
We show that by combining spectral collocation methods, ChC and LGRC, and new efficient numerical methods for solving algebraic MEPs, it is possible to solve such problems both very efficiently and accurately. We improve on several previous results available in the literature, and also present a MATLAB toolbox for solving a wide range of problems.

\subsection{The drift of eigenvalues} \label{drift}

Two techniques are used in order to eliminate the \textquotedblleft{bad}\textquotedblright eigenvalues as well as to
estimate the stability (accuracy) of computations. The first one is the \textit{drift}, with respect to the
order of approximation or the scaling factor, of a set of eigenvalues of interest. In a simplified form this concept has been introduced by J. P. Boyd in \cite{JPB32}. The
second one is based on the check of the \textit{eigenvectors' orthogonality}.

In other words, we want to separate the \textquotedblleft{good}\textquotedblright eigenvalues from the \textquotedblleft{bad}\textquotedblright ones, i.e., inaccurate eigenvalues. An
obvious way to achieve this goal is to compare the eigenvalues computed for different orders of some parameters such as the approximation order (cut-off parameter) $N$ or the scaling factor. Only
those whose difference or \textquotedblleft{resolution-dependent drift}\textquotedblright is \textquotedblleft{small}\textquotedblright can be believed. Actually, in \cite{JPB32} the so called \textit{absolute (ordinal) drift} with respect to the order of approximation has been introduced.

We extend this definition to the following one. The \textit{absolute (ordinal) drift} of the $jth$ eigenvalue with respect to the parameter $\alpha$ is defined as
\begin{equation}
\delta _{j,absolute,\alpha}:=\left\vert \lambda _{j}^{\left( \alpha_{1}\right) }-\lambda
_{j}^{\left( \alpha_{2}\right) }\right\vert , \quad \alpha_{1} \neq \alpha_{2}, \label{eig_drift}
\end{equation}%
where $\lambda _{j}^{\left( \alpha \right) }$ is the $jth$ eigenvalue, after the
eigenvalues have been sorted, as computed using a specific value of the parameter.  In the most common cases this parameter can be $N$ or $c.$

The dependence of $\delta _{j,absolute,\alpha},$ $j=1,2,\ldots, Ne,$ where $Ne$ is the number of analyzed eigenvalues, on the index (mode) $j$ will be displayed in a log-linear plot. If we divide the right hand side of (\ref{eig_drift}) by $\left\vert \lambda _{j}^{\left( \alpha_{1}\right) }\right\vert$ we get the so called \textit{relative drift} denote by $\delta _{j,relative,\alpha}.$

\section{Numerical benchmark problems and discussions} \label{numerics}

\subsection{A regular Schr\"odinger eigenproblem} \label{Regular_S}
The bounded Coffey-Evans potential reads
\begin{equation}
q\left( x\right) :=-2\beta \cos \left( 2x\right) +\beta ^{2}\sin ^{2}\left(
2x\right) ,\ \beta \in \RR.
\label{C-E_pot}
\end{equation}
We attach to equation (\ref{S_eq}) the homogeneous Dirichlet boundary
conditions $u\left( \pm \pi /2\right) =0$ and use $\beta :=30.$
In spite of being regular, the Coffey-Evans
problem is one of the most difficult test problems in the literature because there are very close eigenvalue triplets as $\beta$ increases.

We have used a short Chebfun code in order to solve this regular
eigenproblem. It is available in the next lines. With appropriate changes, this code can be used to analyze any other
Schr\"odinger eigenproblem.
\begin{verbatim}
dom=[-pi/2,pi/2];
x=chebfun('x',dom);beta=30; sigma=-1;
L=chebop(dom);
L.op =@(x,y) -diff(y,2)+(-2*beta*cos(2*x)+(beta*sin(2*x))^2)*y;
L.rbc=0; L.lbc=0; N=201;
[V,D]=eigs(L,N,sigma); D=diag(D)
\end{verbatim}

The eigenvalues obtained by ChC and Chebfun are extremely close, practically indistinguishable. They also compare very well with those computed in \cite{Ledoux1} by some coefficient approximation methods of orders $2$, $4$ and $8$ and reported in Table 3 of this paper.

\begin{figure}
\centering
  \includegraphics[scale=0.75]{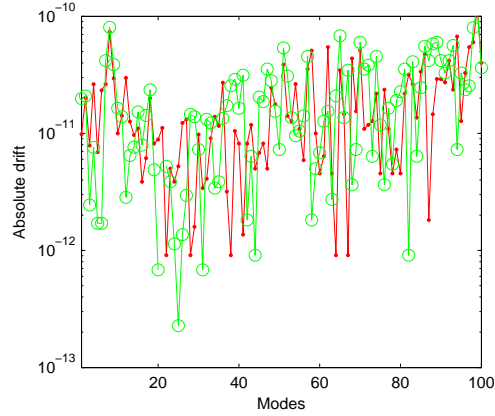}
\caption{The drift of the first $100$ eigenvalues of Coffey-Evans problem. ChC used the orders of approximations $N_{1}:=256$ and $N_{2}:=512$-red dotted line and respectively $N_{1}:=400$ and $N_{2}:=512$-green circled line.}
\label{fig:1}       
\end{figure}

\begin{figure}
\centering
  \includegraphics[scale=0.55]{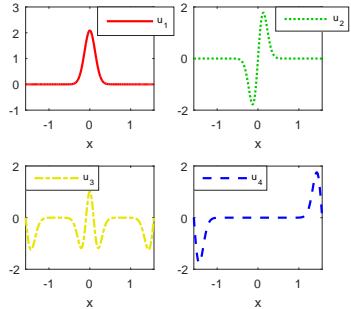}
\caption{The first four eigenvectors of Coffey-Evans problem computed by Chebfun. }
\label{fig:8}       
\end{figure}

\begin{figure}
\centering
  \includegraphics[scale=0.85]{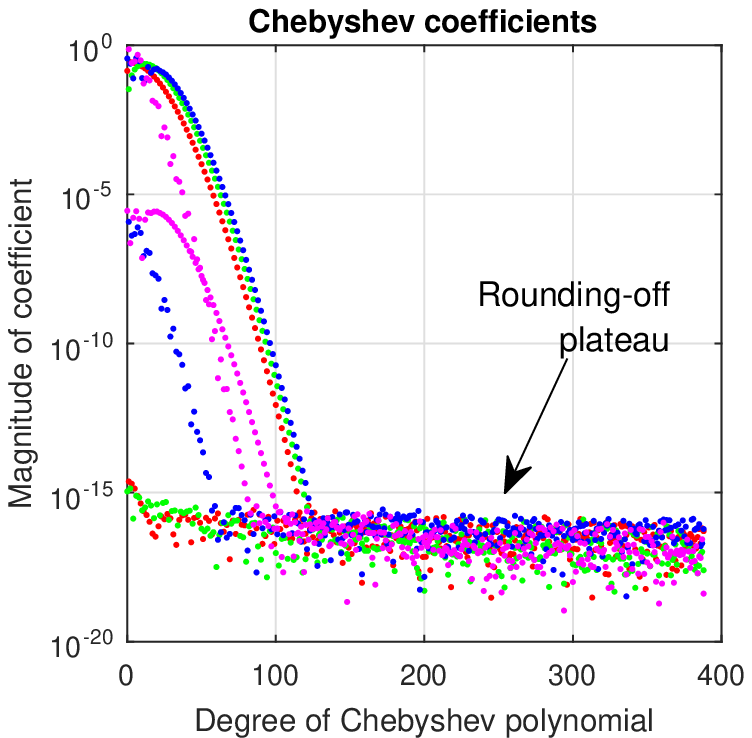}
\caption{The Chebyshev coefficients of the first four eigenvectors of Coffey-Evans problem computed by Chebfun.}
\label{fig:9}       
\end{figure}

The absolute drift reported in Fig. \ref{fig:1} means that we can compute the first hundred eigenvalues with better accuracy than $10^{-10}$. Unfortunately no  accuracy analysis is reported in \cite{Ledoux1} (see also \cite{Ledoux} and \cite{Ledoux06}).
The eigenvectors are either symmetric (the even ones) or
anti-symmetric (the odd ones). The first four of them are depicted in Fig. \ref{fig:8}. They look fairly smooth and satisfy the boundary conditions.
The Chebyshev coefficients of the first four eigenvectors of Coffey-Evans problem computed by Chebfun are displayed in Fig. \ref{fig:9}. These coefficients decrease sharply and smoothly to a rounding-off plateau below $10^{-15}$. Roughly speaking this means they are computed with the machine precision.
\begin{table}[htbp]
\centering
\caption{High index eigenvalues of Schr\"odinger eigenproblem equipped with Coffey-Evans potential (\ref{C-E_pot}) computed by three different methods.}
\label{tab:2}
\begin{tabular}{cccc}
\hline
$\;j$ & $\;\lambda_{j}$ by Chebfun & $\;\lambda_{j}$ computed in \cite{Ledoux1} & $\;\lambda_{j}$ by ChC  \\\hline
$\;20$ & $\; \numprint{951.8788067958783}$ & $\; 951.878806796591$ & $\; 951.8788067965993$\\
$\;30$ & $\; 1438.295244640637$ & $\; 1438.295244640802$ & $\; 1438.295244640797$\\
$\;40$ & $\; 2146.405360539156$ & $\; 2146.405360539854$ & $\; 2146.405360539845$\\
$\;50$ & $\; 3060.923491511540$ & $\; 3060.923491511421$ & $\; 3060.923491511401$\\
$\;100$ & $\;10653.52543568510$ & $\; 10653.525435875921$ & $\; 10653.52543587600$\\
$\;200$ & $\;40851.63764596094$ & $\; 40851.637646050455$ & $\; 40851.63764605047$\\\hline
\end{tabular}%
 \end{table}
For this potential
the first eigenvalue  $\lambda_{0}$ is close to zero (actually we have got $\lambda_{0}=\numprint{-6.254959429708980e-12}$) and there are very close eigenvalue triplets
$\left( \lambda_{2},\; \lambda_{3},\; \lambda_{4} \right)$, $\left( \lambda_{6},\; \lambda_{7},\; \lambda_{8} \right)$, $\ldots$ as $\beta$ increases. The common numeric part of the eigenvalues in the first triplet is $2.31664929$ and that of the eigenvalues in the second is a little shorter, i.e. $4.45283.$

\subsection{Two singular Schr\"odinger eigenproblem on the half line} \label{Singular_S_half}

In this section we study the \textit{mixed spectrum} of some Schr\"odinger eigenproblems having a \textquotedblleft potential well\textquotedblright dying out at infinity.

\subsubsection{Hydrogen atom equation}
The first example consists in the equation (\ref{S_eq}) equipped with the potential
\begin{equation}
q\left( x\right) :=-\frac{1}{x}+\frac{l\left( l+1\right) }{x^{2}},\ l\in \RR,\label{H_pot}
\end{equation}
along with the boundary conditions
\begin{equation}
u(0)=0,\ u\rightarrow 0\ as\ x\rightarrow \infty .  \label{H_bc}
\end{equation}
The problem (\ref{S_eq})-(\ref{H_pot})-(\ref{H_bc}) is clearly singular. We must mention from the beginning that our numerical experiments performed with LGRC and Chebfun together with the truncation of the domain did not produce satisfactory results. In these conditions we have resorted to the \textit{mapped ChC method}.

In order to implement this method we use the \textit{algebraic map}
\begin{equation}
x:=c\frac{1+s}{1-s},\ s\in \left[ -1,1\right] ,\ x\in \lbrack
0,+\infty ), c \in \RR, c > 0, \label{transf}
\end{equation}
which, for each $c$, transforms the interval $\left[ -1,1\right]$ into the half line, and its inverse. The parameter $c$ is free to be tuned
for optimum accuracy.

The mapping (\ref{transf}) has been introduced in \cite{Schonfelder} where its practical effects have been discussed.
The author observed that the convergence of the Chebyshev expansion is governed by the closeness of the singularities of the function being expanded to the expansion region. The major effect of such mapping is to allow
us to move the singularities further away from the expansion region. The mapping
may also weaken the effect of the singularities by modifying the strength of the singularity
as well as moving it. The value of the scaling parameter $c$ used in our computation was chosen essentially by trial and
error along with the drift with respect to this parameter.

In order to write down any second order differential
equation in independent variable $s$, we need the following derivatives:
\begin{equation}
\begin{array}{c}
u^{\prime }\left( x\right) =u^{\prime }\left( s\right)
\frac{1}{x'_{s}},\
x'_{s}=\frac{2c}{\left( 1-s\right) ^{2}}, \\
u^{\prime \prime }\left( x\right) =u^{\prime \prime }\left( s\right) \frac{1%
}{\left( x'_{s}\right) ^{2}}-u^{\prime }\left( s\right)
\frac{x_{s}^{\prime
\prime }}{\left( x'_{x}\right) ^{3}}.
\end{array}
\label{trans_deriv}%
\end{equation}
Now it is easy to write the differential equation for $u\left(s\right)$ and to attach to this new equation the homogeneous Dirichlet boundary conditions $u\left(\pm 1\right)=0.$ As usual we implement these conditions by deleting the first and the last rows and columns of the collocation matrix attached to the left hand side of the equation (\ref{S_eq})-(\ref{H_pot}). This is the most simplified version of removing technique of the independent boundary conditions introduced in \cite{JH}.
With (\ref{trans_deriv}) the following MATLAB code solve this problem.
\begin{verbatim}
% approximation order and parameter l
N=1600; l=1;
% number of displayed eigenvalues
Ne=50;
c=2;         % scaling factor
% Chebyshev differentiation matrices (Weideman & Reddy)
[st,D]=chebdif(N,2);
% N-2 kept nodes; 1 and N are removed nodes
k=2:N-1;  s=st(k);
% enforced boundary conditions in differentiation matrices
D2=D(k,k,2); D1=D(k,k,1);
% collocation matrix of the system
A=-diag((1-x).^4)*D2/(4*(c^2))+diag((1-x).^3)*D1/(2*(c^2))+...
    diag(l*(l+1)*((1-x).^2)./(((1+x).^2)*(c*2))-(1-x)./(c*(1+x)));
% computed and sorted eigenpairs
[U,S]=eig(A); S=diag(S); [t,o]=sort(S); S=S(o); U=U(:,o);
disp(S(1:Ne))
% Chebyshev coefficients of the first four eigenvectors by FCT
Ucoeff=fcgltran(U(:,1:4),1);
\end{verbatim}
\begin{figure}
  \includegraphics{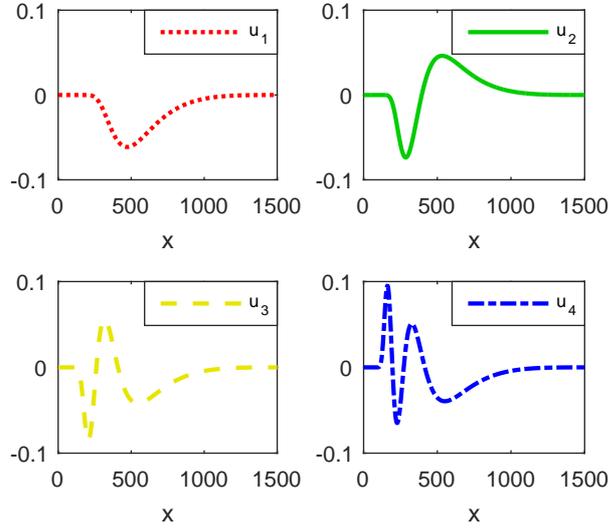}
\caption{The first four eigenvectors for hydrogen atom eigenproblem (\ref{S_eq})-(\ref{H_pot})-(\ref{H_bc}) with $l:=1$ computed by mapped ChC with scaling factor $c:=2$ and $N:=512.$}
\label{fig:3}       
\end{figure}
The first four vectors of the problem are displayed in Fig. \ref{fig:3}. It is clear that they satisfy both boundary conditions but in the right neighborhood of origin they have a totally different behavior from the eigenvectors of regular problems, i.e., they vanish out on continuous portions and not in discrete points. However, they clearly approximate square-integrable eigenfunctions and thus confirm some theoretical results proved in \cite{BR}.
\begin{figure}
\centering
  \includegraphics{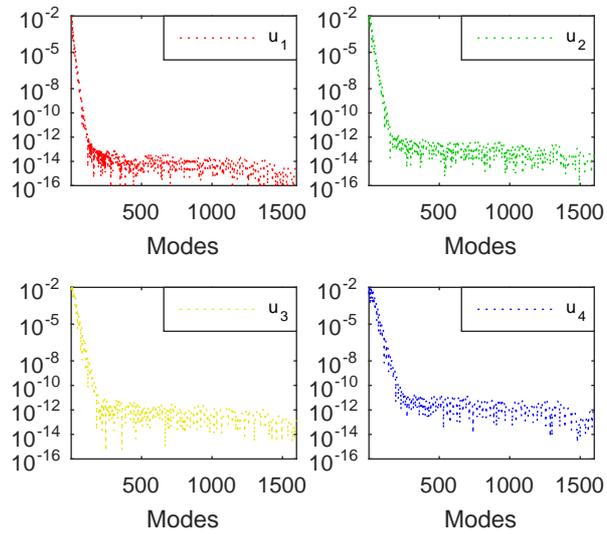}
\caption{The coefficients of the first four eigenvectors for hydrogen atom eigenproblem (\ref{S_eq})-(\ref{H_pot})-(\ref{H_bc}) with $l:=1$ computed by mapped ChC with scaling factor $c:=2$ and $N:=1600.$}
\label{fig:4}       
\end{figure}
Their Chebyshev coefficients obtained using FCT (fast Chebyshev transform-see \cite{GvW} for details) are displayed in Fig.\ref{fig:4}. They decrease sharply and smoothly to some limits, followed by a wide rounding-off plateau. For the first vector (the rightmost one) this limit is around $10^{-13}$. It increases with the index of the vector.
\begin{figure}
\centering
  \includegraphics{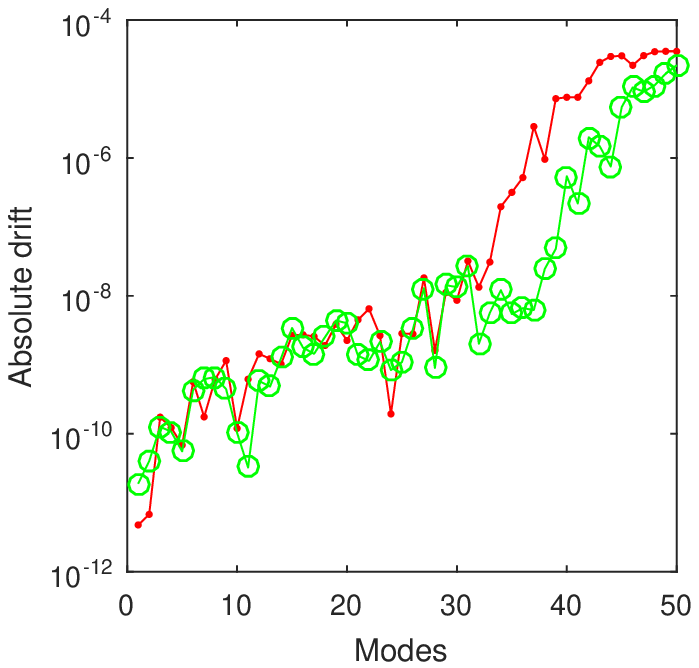}
\caption{The absolute drift with respect to $N$ of the first $50$ negative eigenvalues of hydrogen atom eigenproblem (\ref{S_eq})-(\ref{H_pot})-(\ref{H_bc}) computed by mapped ChC with scaling factor $c:=2$ and $N_{1}:=1600$ and $N_{1}:=2048$ (red line). The green circled line signifies the drift of eigenvalues computed by mapped ChC using $N:=2048$ with respect to the exact eigenvalues $\lambda_{n}=-1/\left(2n+4\right)^2,$
$n=0,1,2,\ldots,50$ corresponding to $l:=1.$}
\label{fig:5}       
\end{figure}
Roughly this means that we cannot hope for a better approximation than something of the order $10^{-13}$ when computing the first eigenvector.

When $N:=2048$ mapped ChC has found $\lambda_{0}=\numprint{-6.250000000166379e-02}$ which is a very good approximation of $-1/16.$ The largest negative eigenvalue has been $\lambda_{66}=\numprint{-7.865782521027431e-07}$ and the next eigenvalue, the smallest positive has been computed as $\lambda_{67}=\numprint{6.517056834998433e-05}.$
To see how the mapped ChC method simulates the notion of continuous spectrum we will provide in Table \ref{tab:3} some significant eigenvalues.
\begin{table}[htbp]
\centering
\caption{High index eigenvalues of hydrogen atom eigenproblem (\ref{S_eq})-(\ref{H_pot})-(\ref{H_bc}).}
\label{tab:3}
\begin{tabular}{cc}
\hline
$\;j$ & $\;\lambda_{j}$ by mapped ChC   \\
\hline
$\;70$ & $\; \numprint{3.346710010488799e-4}$ \\
$\;80$ & $\; 2.246853675452558e-3$ \\
$\;90$ & $\; 6.551029735506819e-3$ \\
$\;100$ & $\; 1.474442763764248e-2$ \\
$\;110$ & $\; 2.890880952550005e-2$ \\
$\;120$ & $\; 5.183912194697291e-2$ \\
\hline
\end{tabular}%
 \end{table}

In the monograph \cite{BR} Sect. 18 it is proved that if the potential $q\left(x\right)$ is continuous and has the asymptotic behaviour $q\left(x\right)=\frac{A}{x}+O\left(\frac{1}{x^2}\right)$ as $x\rightarrow \infty,$ for $\lambda>0$ the spectrum is continuous and the eigenfunctions are not square-integrable and for $\lambda<0$ the spectrum is discrete and the eigenfunctions are square-integrable. The numerical results gathered around this problem plainly confirm this analytical result.

\begin{figure}[htbp]
\centering
  \includegraphics{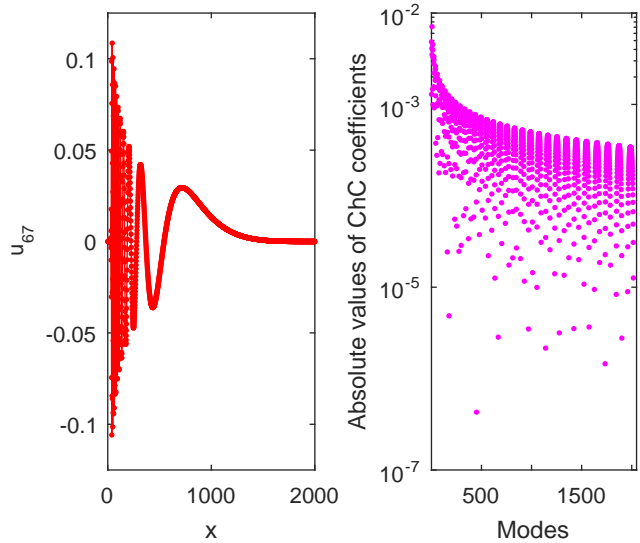}
\caption{(\textbf{a}) The eigenvector corresponding to the first positive eigenvalue of hydrogen atom eigenproblem (\ref{S_eq})-(\ref{H_pot})-(\ref{H_bc})  computed by mapped ChC with scaling factor $c:=2$ and $N:=2024$, $l:=1.$ (left panel). (\textbf{b})In the right panel, in a log-linear plot we display the Chebyshev coefficients of this eigenvector. }
\label{fig:10}       
\end{figure}

Actually the eigenvector corresponding to the first positive eigenvalue (the sixty-seventh one) is depicted in the left panel of Fig. \ref{fig:10}. It is hard to believe that this could approximate a square-integrable eigenfunction. In the right panel of the same figure we displayed the Chebyshev coefficients of the corresponding eigenvector. The oscillations of these coefficients mimic the steep oscillations of the eigenvector.

For the real \textit{scaling factor $c$} we have to mention that in case of eigenvalue problems it can be adjusted only on the mathematical basis. Thus, it can be tuned in order to:
\begin{itemize}
\item improve the decaying rate of the coefficients of spectral expansions.
\item find as orthogonal as possible eigenvectors to an eigenproblem;
\end{itemize}

\subsubsection{Potential with a Coulomb type decay}
\begin{figure}[htbp]
  \includegraphics[scale=0.85]{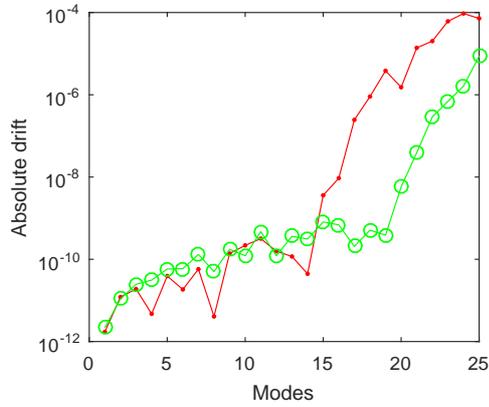}
\caption{The absolute drift with respect to $N$ of the first $25$ eigenvalues for the Schr\"odinger problem (\ref{S_eq})-(\ref{H_bc}) with Coulomb-type decay potential
(\ref{C_pot}), computed by mapped ChC with scaling factor $c:=2$; red dotted line for the case $N_{1}:=512$ and $N_{2}:=1024$ and green circled line for the case $N_{1}:=760$ and $N_{2}:=1024.$}
\label{fig:6}
\end{figure}

Let's consider now the Schr\"odinger equation (\ref{S_eq}) equipped with the potential (\ref{C_pot})
\begin{equation}
q\left( x\right) :=-\frac{1-5\exp \left( -2x\right) }{x}+\frac{l\left(
l+1\right) }{x^{2}},\ l\in \RR.
\label{C_pot}
\end{equation}%
and supplied with
boundary conditions (\ref{H_bc}).
Some eigenvalues computed by ChC method when $c:=2$ and $N:=512$ are compared in Table \ref{tab:1} with their counterparts computed by LGRC and perturbation methods. If for the first eigenvalues the coincidence is excellent, the same does not happen for higher indices. That is why in Fig. \ref{fig:6} in a log-linear plot we illustrate the absolute drift in the case of the ChC method. It is observed that, for instance for the tenth mode, the difference between the values of this eigenvalue calculated with two very different orders of approximation N is of the order $10^{-10}$. This leads us to believe that the eigenvalues computed with ChC are the most accurate. It is also clear from this figure that in the second case, i.e., $N_{1}:=760$ and $N_{2}:=1024$ the drift oscillates less than in the first case reported in this figure.

\begin{table}[htbp]
\centering%
\caption{The first five eigenvalues and the tenth one of the Schr\"odinger problem (\ref{S_eq})-(\ref{H_bc})  when the potential
has a Coulomb-type decay (\ref{C_pot}), computed by three different methods.}
 \label{tab:1}
\begin{tabular}{cccc}
\hline
$\;j$ & $\;\lambda_{j}$ computed by LGRC in \cite{Cig18} & $\;\lambda_{j}$ computed in \cite{Ledoux06} & $\;\lambda_{j}$ by mapped ChC  \\
\hline
$\;1$ & $\; \numprint{-0.0616818466333}$ & $\; -0.061681846633$ & $\; -0.06168184663316705$\\
$\;2$ & $\; -0.0274980999429$ & $\; -0.027498099943$ & $\; -0.02749809994382280$\\
$\;3$ & $\; -0.0155015616910$ & $\; -0.015501561691$ & $\; -0.01550156169420540$\\
$\;4$ & $\; -0.0099354968508$ & $\; -0.009935496851$ & $\; -0.009935496853885005$\\
$\;5$ & $\; -0.0069067013822$ & $\; -0.006906701382$ & $\; -0.006906701375461869$\\
$\;10$ & $\; -0.001963685230$ & $\; -0.001736111111$ & $\; -0.002059879612641054$\\
\hline
\end{tabular}%
 \end{table}

 \subsection{A singular Schr\"odinger eigenproblem on the real line} \label{Singular_S_line}

Let's consider now the equation (\ref{S_eq}) on the real line with the potential%
\begin{equation}
q\left( x\right) :=x^{2}+\frac{\nu x^{2}}{1+\mu x^{2}},
\label{gen_oscill}
\end{equation}%
where $\mu $ and $\nu $ are real parameters. This is a more general (anharmonic) oscillator than  the simpler harmonic one.

Two \textit{behavioral }%
boundary conditions requiring the boundedness of the solutions at large
distance, i.e., $x\rightarrow \pm \infty $ are attached to this equation.
Actually we impose the conditions
\begin{equation}
u\left( x\right) \rightarrow 0\;  \textrm {as}\;
x\rightarrow \pm \infty .\label{inf_cond}
\end{equation}

We have to observe that these conditions are
automatically satisfied in spectral collocation based on Laguerre, Hermite and sinc functions.

\begin{figure}
\centering
  \includegraphics{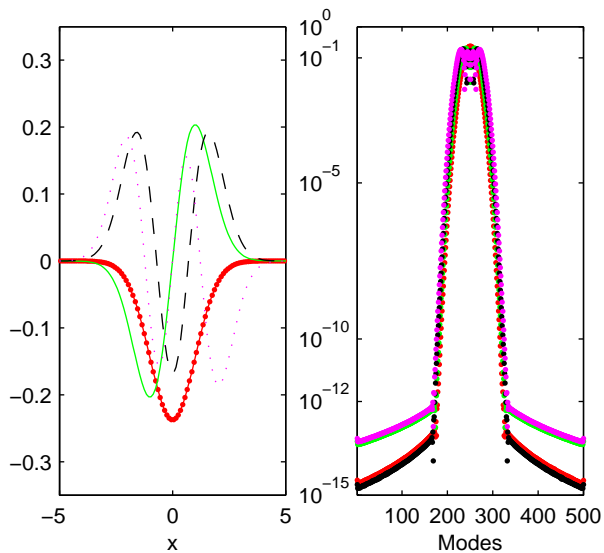}
\caption{(\textbf{a}) Zoom in the first four eigenvectors computed by SiC for the
Schr\"odinger eigenproblem (\ref{S_eq})-(\ref{inf_cond}) with potential (\ref{gen_oscill}) (left panel). (\textbf{b}) In the right panel we display the sinc coefficients for these vectors; red for the first, green for the second, black for the third and magenta for the fourth.}
\label{fig:7}       
\end{figure}

SiC with $N:=500$ and scaling factor $h:=0.1$ has been used in order to produce the following results.
The fist four eigenvectors of Schr\"odinger eigenproblem (\ref{S_eq})-(\ref{inf_cond}) with potential (\ref{gen_oscill}) are
displayed in the left panel of Fig. \ref{fig:7}. Their coefficients are illustrated in the right panel of the same figure. These coefficients symmetrically decrease to approximations of at least $10^{-12}$ which means a reasonable accuracy of the method.
We also have computed some \textit{high-lying} eigenvalues namely $\lambda_{100}$, $\lambda_{150}$ and $\lambda_{200}$, for $\nu:=1$ and $\mu:=500.$ They have respectively the numerical values $$199.001994801512
,\;  301.001995781805,\;  \textrm{and}\;  403.001995224433
.$$ They are five digits approximation for the corresponding harmonic oscillator eigenvalues $\lambda_{n}=2n-1,$ $n=1,2,\ldots .$

\begin{figure}
\centering
  \includegraphics{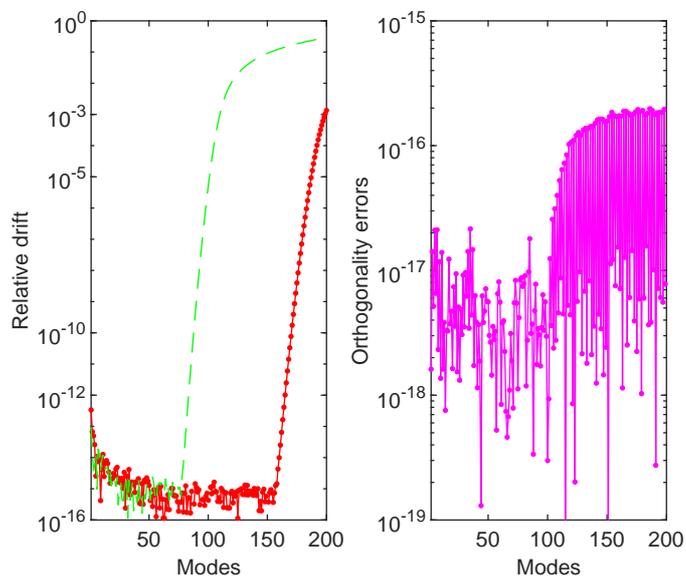}
\caption{(\textbf{a}) The relative drift of the first $200$ eigenvalues of problem (\ref{S_eq})-(\ref{inf_cond}) with potential  (\ref{gen_oscill}), $N_{1}:=500$ and $N_{2}:=400$-red dotted line and $N_{1}:=300$ and $N_{2}:=400$- green line.(left panel). (\textbf{b}) In the right panel we display the orthogonality errors (deficiency) of the first eigenvector of the same problem with respect to the subsequent $199$ of it.}
\label{fig:11}       
\end{figure}

In \cite{Mitra} the author solved this problem numerically for general $\nu$ and $\mu$. He used a very elegant variational argument of Ritz type, based on Hermite functions, and observed that for large $\mu$ and fixed $\nu$ the eigenvalues of this problem approximate those of the harmonic oscillator. We have confirmed this observation even for higher index eigenvalues.
In the left panel of Fig. \ref{fig:11} we display the relative drift of the first $200$ eigenvalues of the problem. It is clear that approximately the first $70$ eigenvalues are calculated with an accuracy better than $10^{-12}$. Eigenvalues with index up to $150$ remain at the same accuracy as $N \approx 500.$ Above this index the accuracy decreases to $10^{-3}$.
Actually, over years, a lot of literature has gathered on this problem. Low order eigenvalues have accurately computed using Runge-Kutta type methods for instance in \cite{Simos} and \cite{Simos1}. In \cite{Trif} the author uses Hermite collocation in order to find only the first eigenvalues for some Schr\"odinger eigenproblem.
 With respect to the departure from orthogonality, i.e. the distance from zero of the scalar product of two eigenvectors, we display in the right hand side of Fig. \ref{fig:11}, in a log linear plot, the absolute values of the of the scalar products of $u_{1}$ and $u_{j}$, $j=2,\ldots,200.$ It is again remarkable that this departure is less than $10^{-15}.$

It is of some importance to justify our choice for SiC. Unlike all other methods of spectral collocation, where the differentiation matrices are highly non-normal, in SiC these matrices are symmetric or skew-symmetric for even respectively odd values of cut-off parameter $N$. And this is an important numerical advantage. For instance in the MATLAB code below we use the routine \textit{eigs} instead of \textit{eig} which would have been more expensive and slower.

The following very simple MATLAB code has been used.
\begin{verbatim}
N=500;                                      % order of approximation
h=0.1;                                      % spacing (scaling factor)
[x,D]=sincdif(N,2,h); D2=D(:,:,2);          % 2nd SiC differentiation
nu=1;mu=500;                                % parameter of the problem
A=-D2+diag((x.^2)+nu*(x.^2)./(1+mu*(x.^2)));% the matrix
[V,D] = eigs(A,250,0); D=diag(D);           % call MATLAB code eigs
[t,o]=sort(real(D)); D=D(o); V=V(:,o);      % sort eigenvalues
\end{verbatim}
This code is fairly similar with that from Section \ref{Singular_S_half} with two differences. First, the Chebyshev differentiation matrices are replaced by sinc differentiation matrices and then the boundary conditions are absent. More exactly, the discrete sinc functions, on which the unknown solution is expanded, are defined by
\begin{equation}
S_{k}\left( x,h\right) :=\frac{\sin \left[ \frac{\pi }{h}\left(
x-x_{k}\right) \right] }{\frac{\pi }{h}\left( x-x_{k}\right) },\
k=1,2,\ldots ,N, \label{sinc_fc}
\end{equation}
and satisfy the boundary conditions (\ref{inf_cond}).
In (\ref{sinc_fc}) the nodes $x_{k}$ are equidistant with spacing $h$ and symmetric with respect to the origin.
 
\section{Concluding remarks and open problems} \label{conclusions}

Regarding the regular problems, Chebfun is unbeatable in terms of accuracy, computation speed, and the information they provide on the accuracy of computational process. It displays the optimal approximation order of unknowns (eigenvectors) and how and to what extent their Chebyshev coefficients decrease. It also specifies the degree to which some boundary conditions are satisfied.

As for the singular problems, the situation is not so offering. In this case, Chebfun leaves enough room for the application of usual (classical) spectral methods. For singular eigenproblems on unbounded domains mapped ChC on the half line and SiC, for problems on the real line, perform better than Chebfun even for Chebfun applied on a mapped domain.

However, approaching these problems in parallel, with Chebfun as well as with the classical spectral methods, one can get greater confidence in the accuracy of numerical results. Regarding ChC method, using the so-called absolute or relative drift in terms of some parameters  we ensure the numerical stability of the numerical process and can eliminate (numerically) spurious eigenvalues. In other words we get automatically and precisely the accuracy at which a specified set of eigenvalues is computed. The departure of orthogonality of a set of eigenvectors provides a useful hint for the accuracy of the numerical process.
We have managed to correctly identify sets of multiple and high indices eigenvalues (triples) and even give some numerical meaning to the notion of continuous spectrum. This second issue obviously remains an open one.

All in all we can say that Chebfun as well as classical spectral methods in various forms, produce more accurate and comprehensive outcomes when solving singular and of course regular eigenproblems than classical methods. However, there remains an open problem, namely that of establishing (possibly automatically) the scaling factor.

The following abbreviations are used in this manuscript:\\

\noindent 
\begin{tabular}{@{}ll}
c & scaling factor for unbounded domains\\
f.e. & finite element method\\
f. d. & finite difference method\\
GPS & generalized pseudospectral method\\
SL & Sturm-Liouville eigenproblem\\
ChC & Chebyshev collocation\\
LGRC & Laguerre Gauss Radau Collocation\\
SiC & sinc collocation\\
MATSLISE & Matlab package for SL and Schr\"odinger equations\\
MEP & Multiparameter Eigenvalue Problem\\
N & the order of approximation of spectral method (cutting-off parameter)\\
SLEDGE & SL Estimates Determined by Global Errors\\
SLEIGN & FORTRAN package for numerical solution to SL eigenproblem\\
SLDRIVER & interactive package for the previous packages\\
\end{tabular}




\end{document}